\documentclass[a4paper,10pt]{article}
\usepackage{authblk}

\usepackage[ngerman,english]{babel}
\usepackage[latin1]{inputenc}
\usepackage{csquotes}
\usepackage{amsfonts,amsmath,amsthm}
\usepackage{empheq}
\usepackage[titletoc,title]{appendix}
\usepackage[backend=bibtex8,doi=false,eprint=false,firstinits=true,isbn=false,style=numeric-comp,maxnames=99]{biblatex}
\makeatletter
\def\blx@maxline{77}
\makeatother
\bibliography{bibl.bib}
\AtEveryBibitem{\clearlist{language}}

\usepackage{cases}
\usepackage{mathabx}
\usepackage{bbm}
\usepackage{xfrac}
\usepackage{fancyhdr}
\usepackage{color}
\usepackage[colorlinks]{hyperref}
\definecolor{blue75}{rgb}{0,0,.75}
\definecolor{green75}{rgb}{0,.75,0}
\hypersetup{colorlinks=true, urlcolor=blue75,linkcolor=blue75,citecolor=green75,pdfstartview=FitB,bookmarksopen=true,bookmarksopenlevel=1}
\usepackage[a4paper, left=2.5cm, right=2.5cm, top=2.5cm,bottom=2cm]{geometry}
\usepackage{constants}
\newcommand{\parenthezises}[1]{\arabic{#1}}
\newconstantfamily{C}{
symbol=C,
format=\parenthezises,
reset={section}
}
\newconstantfamily{M}{
symbol=M,
format=\parenthezises,
reset={section}
}
\usepackage{enumerate}
\usepackage{graphicx}
\graphicspath{ {images/} }
\usepackage{wrapfig}
\usepackage{figbib}
\allowdisplaybreaks
\begin{document}

\newcommand{\cb}{\color{blue}}
\newcommand{\cred}{\color{red}}
\newcommand{\cmg}{\color{magenta}}
\newcommand{\D}{\mathbb{D}}
\newcommand{\R}{\mathbb{R}}
\newcommand{\N}{\mathbb{N}}

\def\diam{\operatorname{diam}}
\def\dist{\operatorname{dist}}
\def\diver{\operatorname{div}}
\def\ess{\operatorname{ess}}
\def\inner{\operatorname{int}}
\def\osc{\operatorname{osc}}
\def\sign{\operatorname{sign}}
\def\supp{\operatorname{supp}}
\newcommand{\BMO}{BMO(\Omega)}
\newcommand{\LOne}{L^{1}(\Omega)}
\newcommand{\LOnen}{(L^{1}(\Omega))^d}
\newcommand{\LTwo}{L^{2}(\Omega)}
\newcommand{\Lq}{L^{q}(\Omega)}
\newcommand{\Lp}{L^{2}(\Omega)}
\newcommand{\Lpn}{(L^{2}(\Omega))^d}
\newcommand{\LInf}{L^{\infty}(\Omega)}
\newcommand{\HOneO}{H^{1,0}(\Omega)}
\newcommand{\HTwoO}{H^{2,0}(\Omega)}
\newcommand{\HOne}{H^{1}(\Omega)}
\newcommand{\HTwo}{H^{2}(\Omega)}
\newcommand{\HmOne}{H^{-1}(\Omega)}
\newcommand{\HmTwo}{H^{-2}(\Omega)}

\newcommand{\LlogL}{L\log L(\Omega)}

\def\avint{\mathop{\,\rlap{-}\!\!\int}\nolimits}

\newtheorem{Theorem}{Theorem}[section]
\newtheorem{Assumption}[Theorem]{Assumptions}
\newtheorem{Corollary}[Theorem]{Corollary}
\newtheorem{Convention}[Theorem]{Convention}
\newtheorem{Definition}[Theorem]{Definition}
\newtheorem{Example}[Theorem]{Example}
\newtheorem{Lemma}[Theorem]{Lemma}
\newtheorem{Notation}[Theorem]{Notation}
\theoremstyle{definition}
\newtheorem{Remark}[Theorem]{Remark}

\newtheorem{proofpart}{Step}
\makeatletter
\@addtoreset{proofpart}{Theorem}
\makeatother
\numberwithin{equation}{section}
\title{Generalised  supersolutions with  mass control for the Keller-Segel  system with logarithmic  sensitivity}
\author{Anna~Zhigun}
\renewcommand\Affilfont{\itshape\small}
\affil{Technische Universität Kaiserslautern, Felix-Klein-Zentrum für Mathematik\\ Paul-Ehrlich-Str. 31, 67663 Kaiserslautern, Germany\\
  e-mail: {zhigun@mathematik.uni-kl.de}}
\date{}
\maketitle
\begin{abstract}
 The existence of  generalised global supersolutions with a control upon the total muss is established for the  parabolic-parabolic Keller-Segel system with logarithmic sensitivity for any space dimension. It is verified that  smooth supersolutions of this sort are actually  classical solutions. 
 Unlike the   previously existing constructions, neither is the  chemotactic sensitivity coefficient required to be small, nor is it necessary for the initial data to be radially symmetric. 
\\\\
{\bf Keywords}: chemotaxis;  generalised supersolution; global existence;  logarithmic sensitivity.\\
MSC 2010: 
35B45, 
92C17,  
35D30, 
35D99, 
35K55. 
\end{abstract}

\section{Introduction}
Coupled reaction-diffusion-transport PDEs are a standard tool in the mathematical modelling of cell motility on the macroscale. Thereby, the diffusion-dominated systems are among the best  studied analytically. Standard theory (see, e.g., \cite{LSU}) ensures the existence of bounded solutions to such systems.  Still, they are not always the optimal choice. Indeed, in many instances it is not the chaotic movement but, rather, the active drift of the cells towards some substance  which actually dominates the motion and, as a result, may lead to a strong aggregation of the biomass density. The best known model example of such a situation is provided by the celebrated Keller-Segel system for chemotaxis \cite{KS70,KS71}. This parabolic-parabolic system and its parabolic-elliptic simplifications have been objects of extensive studies in recent decades.
 It turned out that in higher spatial dimensions the solutions to such systems can exhibit a blow-up in finite time which calls into question the global solvability.  
For certain parabolic-elliptic versions of the classical Keller-Segel model on the plane one was able to extend the solutions which collapse into a persistent Dirac-type singularity beyond a finite-time blow-up by constructing measure-valued continuations  \cite{LuckSugVel,WinklerTello2013}.
For a detailed overview of available results concerning boundedness/blow-up, as well as other properties, of the Keller-Segel model the reader is referred to 
\cite{Horstmann,BBTW}.

Very recently a new  solution concept was introduced  \cite{LankWink2017} in the context of a version of the Keller-Segel system with a signal-dependent chemotactic sensitivity:
\begin{subequations}\label{KSLog}
 \begin{empheq}[left={ \empheqlbrace\ }]{align}
  &\partial_tu=\nabla\cdot\left(\nabla u-\chi \frac{u}{v}\nabla v\right)\quad&&\text{ in }\R^+\times\Omega,\label{Equ}\\
  &\partial_t v=\Delta v-v+u\quad&&\text{ in }\R^+\times\Omega,\label{Eqv}\\
  &\partial_{\nu} u=\partial_{\nu} v=0\quad&&\text{ in }\R^+\times\partial\Omega,\label{bc}\\
  &u(0,\cdot)=u_0,\ v(0,\cdot)=v_0\quad&&\text{ in }\Omega,
 \end{empheq}
\end{subequations}
where $\Omega$ is a smooth bounded domain in $\R^n$, $n\in\N$, with the corresponding outer normal unit vector $\nu$ on $\partial\Omega$, and  $\chi$ is a positive number. In this model prototype the  cells are assumed to respond to the changes of the logarithm of the signal concentration thus  following the Weber-Fechner law. Due to the saturation effect upon the chemotactic sensitivity in  the presence of high levels of signal concentration, the solutions of both \eqref{KSLog} and the corresponding parabolic-elliptic versions are less prone to the formation of strong singularities, such as, e.g., Dirac measures, than  those of the classical Keller-Segel model. In particular, the global existence of bounded classical \cite{FujKent2015,ZhaoZheng2016,Wink2011,MizYokota2017,Lankeit2017,BBTW,FujieSenba2016II,FujieSenba2016,NagaiSenba1998}, weak \cite{StinWink2011,Wink2011}, and   generalised \cite{LankWink2017,Black2018} solutions was established for certain ranges of parameter $\chi$ which depend upon $n$ and, also, on  whether the setting is radial-symmetric or not. On the other hand,  it is known, for a  parabolic-elliptic case at least \cite{NagaiSenba1998} that blow-up solutions exist for $n\geq3$ and $\chi$ large enough. To the best of our knowledge, no solution concept has as as yet been proposed for \eqref{KSLog} which would allow to treat the general non-radial symmetric setting and arbitrary large  $\chi>0$ if $n\geq3$. The present work aims to fill this gap.  

 In \cite{LankWink2017}, the  authors introduced the coupled quantity
\begin{align}
 F(u,v)=u^pv^q\label{Fpos}
\end{align}
for $$p,q\in (0,1)$$ and derived a variational inequality which it should satisfy provided that $u$ and $v$ solve \eqref{KSLog} in the classical sense. They showed that such inequality, when   complemented by a control upon the total mass and equation \eqref{Eqv} fulfilled in the usual  weak sense, together  comprise a reasonable concept of a global  solution to the full parabolic-parabolic system \eqref{KSLog}. This construction made it possible to  extend the range of $\chi$'s for which \eqref{KSLog} can be considered as globally solvable, namely:
\begin{align}
 0<\chi<\begin{cases}
       \infty&\text{if }n=2,\\
       \sqrt{8}&\text{if }n=3,\\
       \frac{n}{n-2}&\text{if }n\geq4
      \end{cases}\qquad \text{\cite{LankWink2017}}.\label{chi}
\end{align}
Indeed, previously available   results on  global  solvability  presupposed either a more restrictive condition on $\chi$ \cite{Wink2011}, or the radial symmetry requirement \cite{StinWink2011}. 
The generalised solution from \cite{LankWink2017} is a limit of a regularising sequence. 
Thanks to condition \eqref{chi}, the $u$-component of this sequence is  uniformly integrable over arbitrary finite time-space cylinders. An important consequence of the later  is  the fact  that the generalised solution satisfies equation \eqref{Eqv} in the usual weak sense. In this paper we further develop  the framework presented in \cite{LankWink2017}.
The key idea is to replace \eqref{Fpos} with 
 \begin{align}
 F(u,v)=u^{-a}v^{-b}\label{F}
\end{align}
where parameters $a$ and $b$ are assumed to satisfy
\begin{align}
 a>0\qquad\text{and}\qquad b>b_+(a):=\frac{1+a}{2}\left(\sqrt{1+\chi^2a}-1\right),\label{ab}%
\end{align}
so that clearly
\begin{align*}
 b>0
\end{align*}
as well. An advantage of switching to a negative power of $u$ in \eqref{F} is that the imposed  restriction   on $b$ in \eqref{ab} can be met for any $\chi>0$. 
On the negative side is that the $u$-components of the approximation sequence  need no longer be locally in time uniformly integrable. As a result, the limit pair $(u,v)$, while remaining a  supersolution to \eqref{Eqv},  may, however, fail to be its subsolution, even in the weak sense.  Thus, our approach to  \eqref{Eqv} is to replace it by an identity which contains a certain nonnegative Radon measure. Notwithstanding, thanks to a suitable mass control involving this measure and an accurate description of the boundary conditions based on the  concept of normal traces for divergence-measure fields \cite{ChenFrid2001}, our generalised supersolution coincides with a classical solution in the case of regular $u$ and $v$. 

The rest of the paper is organised in the following way. In {\it Section \ref{SecSol}} we introduce   and discuss the assumptions on the initial data, the proposed supersolution concept (cl. {\it Definition \ref{defsol}}), and   two main results: {\it Theorem \ref{TheoEx}} deals with existence of such supersolutions, while {\it Theorem \ref{classol}} establishes a connection to the classical solvability. In {\it Section \ref{secreg}} we prepare some  ingredients which are necessary to our proof of  {\it Theorem \ref{TheoEx}} in {\it Section \ref{SecEx}}. {\it Theorem \ref{classol}} is proved in the closing {\it Section \ref{SecCl}}. Finally, {\it Appendix \ref{AppA}} contains some facts on the divergence-measure fields which we use in this paper. 

\section{Generalised supersolutions to \texorpdfstring{\eqref{KSLog}}{}}\label{SecSol}
We assume throughout that the initial data satisfies the following assumptions:
\begin{align}
 &0<u_0,v_0\in L^1(\Omega)\label{uv0L1}
\end{align}
and
\begin{align}
   u_0^{-a}v_0^{-b}\in L^1(\Omega)\qquad\text{ for some }a,b>0.\label{prL1}
\end{align}
\begin{Remark}[Initial conditions]\label{RemIni}~
 \begin{enumerate}
\item  Conditions  \eqref{uv0L1}-\eqref{prL1} are to be compared to  (1.7) in \cite{LankWink2017}. 
We impose less regularity for $u_0$ and $v_0$, as well as allow, e.g., $v_0$ to touch zero at some points in $\overline{\Omega}$.  Yet, since we deal with negative powers of $u_0$, we cannot consider general $u_0\geq0$ which vanish on sets of non-zero Lebesgue measure.    
\item In our existence result, {\it Theorem \ref{TheoEx}} below, we   require \eqref{ab} to hold in addition. 
Observe that $b_+$ is a strictly increasing function, and one readily sees  that 
\begin{align*}
 b_+(a)\underset{a\rightarrow0}{\rightarrow}0.
\end{align*} 
Hence, given any small $b$, we only need  to satisfy $a\in\left(0,(b_+)^{-1}(b)\right)$, so that $a$ can be chosen arbitrary  small as well. On the whole, this makes it easier to satisfy \eqref{ab}. 
\end{enumerate}
  
\end{Remark}

Motivated by an idea from  \cite{LankWink2017} we introduce the following  concept of a generalised supersolution with a mass control:
\begin{Definition}[Generalised supersolution]\label{defsol} Let $(u_0,v_0)$ satisfy conditions \eqref{uv0L1} and \eqref{prL1}.  We call a pair of measurable functions $(u,v):\R_0^+\times\overline{\Omega}\rightarrow \R^+\times \R^+$ a \underline{generalised supersolution} to system \eqref{KSLog} if 
\begin{enumerate}[(i)]
  \item\label{defi} $u\in L^{\infty}(\R^+;L^1(\Omega))$,\  $v\in L^1_{loc}(\R_0^+;W^{1,1}(\Omega))$,\ $v^{-1}\in L_{loc}^{\infty}(\R^+\times\overline{\Omega})$;
  \item $u^{-\frac{{a}}{2}}v^{-\frac{{b}}{2}}\in L^2_{loc}(\R_0^+;H^1(\Omega))\cap L^{\infty}_{loc}(\R_0^+;L^2(\Omega))$,\  $u^{-\frac{{a}}{2}}v^{-\frac{{b}}{2}-1}\nabla v\in L^2_{loc}(\R_0^+;L^2(\Omega))$,\\ $u^{-a+1}v^{{-b}-1}\in L^1_{loc}(\R_0^+;L^1(\Omega))$;
  \item\label{inDMp} $\int_0^{\infty}\psi\left(\nabla\left(u^{-a}v^{-b}\right)+\chi au^{-a}v^{-b-1}\nabla v\right)\,ds$,\ \   $\int_0^{\infty}\psi \nabla v\,ds\in {\cal DM}^{p}(\Omega)$   for all $\psi\in C^1_0(\R_0^+)$ for some $p\in\left(1,\frac{n}{n-1}\right)$;
 \item for all $0\leq\varphi\in C^1(\overline{\Omega})$ and $0\leq\psi\in C^1_0(\R_0^+)$ it holds that
\begin{align}
 &-\int_0^{\infty}\partial_t\psi\int_{\Omega}u^{-a}v^{-b}\varphi\,dxdt-\psi(0)\int_{\Omega}u_0^{-a}v_0^{-b}\varphi\,dx\nonumber\\
  \leq  &-4\int_0^{\infty}\psi\int_{\Omega}\left(\frac{a+1}{a}\left|\nabla\left(u^{-\frac{{a}}{2}}v^{-\frac{{b}}{2}} \right)\right|^2+\left(\frac{b}{a}+\chi\frac{a+1}{2}\right)\nabla\left(u^{-\frac{{a}}{2}}v^{-\frac{{b}}{2}} \right)\cdot u^{-\frac{{a}}{2}}v^{-\frac{{b}}{2}-1}\nabla v\right.\nonumber\\
 &\left.\phantom{\int_0^t\psi\int_{\Omega}-4}+\frac{1}{4}\left(\frac{b^2}{a}+b+\chi b\right)\left|u^{-\frac{{a}}{2}}v^{-\frac{{b}}{2}-1}\nabla v\right|^2\right)\varphi\,dxds\nonumber\\
 &-\int_0^{\infty}\psi\int_{\Omega}\left(\nabla\left(u^{-a}v^{-b}\right)+\chi au^{-a}v^{-b-1}\nabla v\right)\cdot\nabla\varphi\,dxds\nonumber\\
 &+\int_0^{\infty}\psi\int_{\Omega}\left(b u^{-a}v^{-b}-b u^{-a+1}v^{{-b}-1}\right)\varphi\,dxds;
 \label{superu}
\end{align}
\item 
there exists  some non-negative Radon measure $\mu$ in $\R_0^+\times \overline{\Omega}$, s.t. for all $0\leq\varphi\in C^1(\overline{\Omega})$ and $0\leq\psi\in C^1_0(\R_0^+)$ it holds that
\begin{align}
 -\int_0^{\infty}\partial_t\psi\int_{\Omega}v\varphi\,dxds-\psi(0)\int_{\Omega}v_0\varphi\,dx=\int_0^{\infty}\psi\int_{\Omega}-\nabla v\cdot\nabla\varphi+(-v+u)\varphi\,dxds+\int_0^{\infty}\psi\int_{\overline{\Omega}}\varphi \,d\mu(s,x),\label{superv}
\end{align}
and 
\begin{align}
 \int_0^{\infty}\psi\|u(t,\cdot)\|_{L^1(\Omega)}\,ds+\int_0^{\infty}\psi\int_{\overline{\Omega}} \,d\mu(s,x)\leq \|u_{0}\|_{L^1(\Omega)}\|\psi\|_{L^1(\R^+)};\label{massu0}
\end{align}
\item
for all $\psi\in C^1_0(\R_0^+)$ it holds that 
\begin{alignat}{3}
&\int_0^{\infty}\psi\left(\nabla\left(u^{-a}v^{-b}\right)+\chi au^{-a}v^{-b-1}\nabla v\right)\,ds\cdot\nu|_{\partial\Omega}=0&&\qquad\text{in }W^{-\frac{1}{p},p}(\partial\Omega),\label{prweakbc}\\
 &\int_0^{\infty}\psi\nabla v\,ds\cdot\nu|_{\partial\Omega}=0&&\qquad\text{in }W^{-\frac{1}{p},p}(\partial\Omega).\label{vweakbc}
\end{alignat}
\end{enumerate}
\end{Definition}
Several remarks are in order. 
\begin{Remark}[Boundary conditions] The variational reformulations \eqref{prweakbc}-\eqref{vweakbc} of the boundary conditions \eqref{bc} are consistent with the regularity assumptions in \eqref{inDMp}, cl. \cite[Theorem 2.1]{ChenFrid2001}. Some necessary facts on the spaces   ${\cal DM}^{p}(\Omega)$ are recalled  for the reader's convenience in  {\it Appendix \ref{AppA}}. 
\end{Remark}
\begin{Remark}[Coupled quantities]~
\begin{enumerate}
\item 
 One of the well-known difficulties  to be faced while dealing with a system like \eqref{KSLog}  in higher dimensions is the  generally poor regularity of  variable $u$. The reason lies in a comparative weakness of  linear diffusion  which often fails to compensate the aggregation due to taxis. 
 At the same time,  variable $v$ is much better-behaved. Our previous studies of  highly-degenerate haptotaxis systems \cite{ZSH,ZSU} suggest that   it can be helpful to introduce a quantity which involves {\it both} variables and possesses an integrable gradient. Thereby, the choice of such coupled quantity  depends to a large extent upon  estimates  one is able to derive for a given system. 
 In particular, in the present case  it turns out to be fruitful to  study  the term $u^{-\frac{{a}}{2}}v^{-\frac{{b}}{2}}$ and to include its  gradient into the  supersolution concept. 
 \item Since $u^{-\frac{{a}}{2}}v^{-\frac{{b}}{2}}\in L^2_{loc}(\R_0^+;H^1(\Omega))$, we have due to the weak chain rule and the H\"older inequality that $u^{-\frac{{a}}{2}}v^{-\frac{{b}}{2}}\in L^1_{loc}(\R_0^+;W^{1,1}(\Omega))$. On the whole, the regularity  imposed by \eqref{defi}-\eqref{inDMp} ensures that all integrals in the variational formulations \eqref{superu}-\eqref{vweakbc} do make sense. 
 \end{enumerate}

\end{Remark}

\begin{Remark}[Comparison with the solution concept from \cite{LankWink2017}]
The main difference between our construction and the generalised solution concept in \cite{LankWink2017} is the presence of a nonnegative  and, in general, nonzero  measure $\mu$ on the right-hand side in the weak formulation \eqref{superv}. This means that the pair   $(u,v)$ is actually a weak   supersolution of \eqref{Eqv} but  may, however, fail to be its subsolution. In particular,
 even if both $u$ and $v$ are smooth functions,  the variational properties \eqref{superu}-\eqref{superv} do not on their own imply that $(u,v)$ is a supersolution to \eqref{superu}. Fortunately,  as it turns out, this property can be saved by taking into account the  the boundary conditions in form of \eqref{prweakbc}-\eqref{vweakbc}, see the proof of {\it Theorem \ref{classol}} below.                                                                                                                                                                                                                                                                                                                                                                                                     

\end{Remark}
Our result on existence now reads:
\begin{Theorem}[Existence of generalised supersolutions]\label{TheoEx}
 Let $\chi$ be any positive number. Let the initial conditions $u_0$ and $v_0$ satisfy \eqref{uv0L1}-\eqref{prL1} for some constants $a$ and $b$ which fulfil  \eqref{ab}. Then there exists a generalised supersolution $(u,v)$ in terms of {\it Definition \ref{defsol}}.
\end{Theorem}
The proof of this theorem is based on a suitable regularisation and a series of priori estimates in {\it Section \ref{secreg}} leading into  a limit procedure in {\it Section \ref{SecEx}}. 

Our interest in the introduced concept of generalised supersolutions is supported by the following result:
\begin{Theorem}[Classical solutions]\label{classol} 
 Let a pair $(u,v)$ be a supersolution in terms of {\it Definition \ref{defsol}}. Assume in addition that 
 \begin{align*}
u,v\in C^{1,2}(\R^+_0\times \overline{\Omega})                                                                                                                                             \end{align*}
 and   
 \begin{align}
\frac{u}{v}\in L^{\infty}_{loc}(\R^+_0\times \overline{\Omega}).\label{udivv}                                                                                                                                                                                                                  \end{align}
 Then $(u,v)$ solves  \eqref{KSLog} in the classical sense. 
\end{Theorem}
The proof of {\it Theorem \ref{classol}} is given in {\it Section \ref{SecCl}}.

\begin{Remark}[Notation]
 We make the following useful convention: For any index $i$, a quantity $C_i$ denotes a positive constant or, alternatively, a positive function, 
which is non-decreasing in each of its arguments.
Moreover, dependence upon such parameters as: the space dimension $n$, domain $\Omega$, constants $a,b,\chi$, as well as the structure of the initial data $u_0,v_0$, is mostly 
{\bf not} indicated in an explicit way.
\end{Remark}

\section{Smooth regularisations for  \texorpdfstring{\eqref{KSLog}}{}}\label{secreg}
Let  
\begin{align*}
2\leq k\in\N.                                             \end{align*}
 Following \cite{LankWink2017}, we consider the a family of regularisations of system \eqref{KSLog}: 
\begin{subequations}\label{KSLoge}
 \begin{empheq}[left={ \empheqlbrace\ }]{align}
  &\partial_t u_{k}=\nabla\cdot\left(\nabla u_{k}-\chi \frac{u_{k}}{v_{k}}\nabla v_{k}\right)\quad&&\text{ in }(0,T)\times\Omega,\label{ue}\\
  &\partial_t v_{k}=\Delta v_{k}-v_{k}+\frac{u_{k}}{1+\frac{1}{k} u_{k}}\quad&&\text{ in }(0,T)\times\Omega,\label{ve}\\
  &\partial_{\nu} u_{k}=\partial_{\nu} v_{k}=0\quad&&\text{ in }(0,T)\times\partial\Omega,\label{uvebc}\\
  &u_{k}(0,\cdot)=u_{k0},\ v_{k}(0,\cdot)=v_{k0}\quad&&\text{ in }\Omega.
 \end{empheq}
\end{subequations}
Thereby, we choose the regularised initial data  $u_{k0}$ and $v_{k0}$ so as to satisfy
\begin{align}
&0<u_{k0},v_{k0}\in W^{1,\infty}(\Omega),\label{RegIni}
\end{align}
as well as to be suitable approximations to the original starting values, $u_0$ and $v_0$, see next sequel. 
Classical theory  for upper-triangular systems (see, e.g., \cite{Amann1}) implies that \eqref{KSLoge} possesses a unique global classical solution  $\left(u_{k},v_{k}\right)$. Moreover, due to the maximum principle,  both solution components are strictly positive in $\R^+_0\times \overline{\Omega}$. These solutions are studied in  {\it Subsections \ref{basic}-\ref{keyest}}.  
\subsection{Approximation of initial data}
Apart from being smooth, we assume  that $u_{k0}$ and $v_{k0}$ fulfil  the following conditions: for all $k\in\N$ it holds that
\begin{alignat}{3}
&u_{k0}\underset{k\rightarrow\infty}{\rightarrow}u_0&&\qquad\text{in }L^1(\Omega)\text{ and a.e. in }\Omega,\label{aproxiniu_}\\
&v_{k0}\underset{k\rightarrow\infty}{\rightarrow}v_0&&\qquad\text{in }L^1(\Omega)\text{ and a.e. in }\Omega,\label{aproxiniv_}\\
&u_{k0}^{-a}v_{k0}^{-b}\underset{k\rightarrow\infty}{\rightarrow}u_0^{-a}v_0^{-b}&&\qquad\text{in }L^1(\Omega)\text{ and a.e. in }\Omega.\label{aproxiniuv_}
\end{alignat}
 Let us check that these conditions can be met. First, we observe that since $u_0,v_0\in L^1(\Omega)$, there exist some sequences $u_{k0}$ and $v_{k0}$ which satisfy \eqref{RegIni} and \eqref{aproxiniu_}-\eqref{aproxiniv_}, as well as
 \begin{align}
  &u_{k0}\geq\begin{cases}
   \max\left\{k^{-\frac{b}{a}},u_0\right\}&\text{ for }u_0<k,\\
         k&\text{ for }u_0\geq k,
         \end{cases}\label{Um0leq}\\
  &v_{k0}\geq\begin{cases}
   \max\left\{k^{-\frac{a}{b}},v_0\right\}&\text{ for }u_0<k,\\
         k&\text{ for }v_0\geq k.
 \end{cases}\label{Vm0leq}
 \end{align}
%
Using \eqref{Um0leq}-\eqref{Vm0leq}, we compute that
\begin{align}
 u_{k0}^{-a}v_{k0}^{-b}\leq& \max\left\{u_0^{-a}v_0^{-b},1\right\}\qquad\text{a.e. in }\Omega.\label{maj}
\end{align}
Combining \eqref{prL1}, \eqref{aproxiniu_}-\eqref{aproxiniv_}, and \eqref{maj} with the dominated convergence theorem we obtain \eqref{aproxiniuv_}. 

\subsection{Basic properties of    \texorpdfstring{\eqref{KSLoge}}{}}\label{basic}
Integrating equations \eqref{ue} and \eqref{ve} over $\Omega$ and using the boundary conditions and partial integration we obtain the following information about the total masses: for all $t\geq0$
\begin{align}
 &\|u_{k}(t,\cdot)\|_{L^1(\Omega)}= \|u_{k 0}\|_{L^1(\Omega)},\label{massue}\\
 &\|v_{k}(t,\cdot)\|_{L^1(\Omega)}\leq\left(1-e^{-t}\right)\|u_{k 0}\|_{L^1(\Omega)}+e^{-t}\|v_{k 0}\|_{L^1(\Omega)}.\label{massve}
\end{align}
Due to \eqref{massue}-\eqref{massve} and a classical result based on duality (see, e.g., the proof of Lemma 5 in \cite[Appendix A]{BOTHE2010120}) we have for all \begin{align}
(r,s)\in\left[1,\frac{n+2}{n}\right)\times\left[1,\frac{n+2}{n+1}\right)\nonumber                                                                                                                                                                \end{align}
that
\begin{align}
 \left\{\left(v_{k},\nabla v_{k}\right)\right\}_{k\in(0,1]} \quad \text{is precompact in }L^r((0,T)\times\Omega)\times L^s((0,T)\times\Omega).\label{vcomp}
\end{align}

Further, using the maximum principle and the strict positivity of the Neumann heat kernel, we conclude from   \eqref{ve} and \eqref{aproxiniv_} that $v_{k}$ can be controlled  from below in the following way: 
\begin{align}
 \inf_{(\tau,T)\times \Omega}v_{k}\geq &\inf_{(\tau,T)\times \Omega}e^{-t}e^{t\Delta}v_{k 0}\nonumber\\
 \geq &\C(\tau,T)\|v_{k 0}\|_{L^1(\Omega)}\nonumber\\
 \geq &\Cl{vmin}(\tau,T)>0\qquad\text{for all }0<\tau<T<\infty.\label{estvmin}
\end{align}
\subsection{A variational formulation for  \texorpdfstring{\eqref{KSLoge}}{}}
Let $F\in C^1(\R^+\times \R^+)$. Multiplying  \eqref{ue} and \eqref{ve} by $\partial_{u}F(u_{k},v_{k})$ and $\partial_{v}F(u_{k},v_{k})$, respectively, adding the results together, and using the chain rule where necessary, we compute  that 
\begin{align}
 &\partial_tF(u_{k},v_{k}) \nonumber\\
 =&-\left(\nabla u_{k}-\chi \frac{u_{k}}{v_{k}}\nabla v_{k}\right)\left(\partial_{uu}F(u_{k},v_{k})\nabla u_{k}+\partial_{uv}F(u_{k},v_{k})\nabla v_{k}\right)
 +\nabla\cdot\left(\partial_{u}F(u_{k},v_{k})\left(\nabla u_{k}-\chi \frac{u_{k}}{v_{k}}\nabla v_{k}\right)\right)\nonumber\\
 &-\nabla v_{k} \cdot\left(\partial_{uv}F(u_{k},v_{k})\nabla u_{k}+\partial_{vv}F(u_{k},v_{k})\nabla v_{k}\right)
 +\nabla\cdot\left(\partial_{v}F(u_{k},v_{k})\nabla v_{k}\right)\nonumber\\
 &+\partial_{v}F(u_{k},v_{k})\left(-v_{k}+\frac{u_{k}}{1+\frac{1}{k} u_{k}}\right)\nonumber\\
 =&-\left(\partial_{uu}F(u_{k},v_{k})|\nabla u_{k}|^2+2\left(\partial_{uv}F(u_{k},v_{k})-\frac{\chi}{2} \frac{u_{k}}{v_{k}}\partial_{uu}F(u_{k},v_{k})\right)\nabla u_{k}\cdot\nabla v_{k}\right.\nonumber\\
 &\phantom{aa}\left.+\left(\partial_{vv}F(u_{k},v_{k})-\chi \frac{u_{k}}{v_{k}}\partial_{uv}F(u_{k},v_{k})\right)|\nabla v_{k}|^2\right)\nonumber\\
 &+\nabla\cdot\left(\partial_{u}F(u_{k},v_{k})\nabla u_{k}+\left(\partial_{v}F(u_{k},v_{k})-\chi \frac{u_{k}}{v_{k}}\partial_{u}F(u_{k},v_{k})\right)\nabla v_{k}\right)+\partial_{v}F(u_{k},v_{k})\left(-v_{k}+\frac{u_{k}}{1+\frac{1}{k} u_{k}}\right).\label{KEYI}
\end{align}
Now we choose $F$ as in \eqref{F}. 
Using the chain rule where necessary, one easily verifies that
\begin{align}
 &\nabla u_{k}=-\frac{2}{a}u_{k}^{\frac{{a}}{2}+1}v_{k}^{\frac{{b}}{2}}\left(\nabla\left(u_{k}^{-\frac{{a}}{2}}v_{k}^{-\frac{{b}}{2}} \right)+\frac{b}{2}\left(u_{k}^{-\frac{{a}}{2}}v_{k}^{-\frac{{b}}{2}-1}\nabla v_{k}\right)\right),\label{nue}\\
 &\nabla v_{k}=u_{k}^{\frac{{a}}{2}}v_{k}^{\frac{{b}}{2}+1}\left(u_{k}^{-\frac{{a}}{2}}v_{k}^{-\frac{{b}}{2}-1}\nabla v_{k}\right).\label{nve}
\end{align}
Plugging \eqref{nue}-\eqref{nve} into \eqref{KEYI} we arrive, after some computation, at the following identity: 
\begin{align}
\partial_t\left(u_{k}^{-a}v_{k}^{-b}\right) 
  = &-4\left(\frac{a+1}{a}\left|\nabla\left(u_{k}^{-\frac{{a}}{2}}v_{k}^{-\frac{{b}}{2}} \right)\right|^2+\left(\frac{b}{a}+\chi\frac{a+1}{2}\right)\nabla\left(u_{k}^{-\frac{{a}}{2}}v_{k}^{-\frac{{b}}{2}} \right)\cdot u_{k}^{-\frac{{a}}{2}}v_{k}^{-\frac{{b}}{2}-1}\nabla v_{k}\right.\nonumber\\
 &\phantom{-4aa}\left.+\frac{1}{4}\left(\frac{b^2}{a}+b+\chi b\right)\left|u_{k}^{-\frac{{a}}{2}}v_{k}^{-\frac{{b}}{2}-1}\nabla v_{k}\right|^2\right)\nonumber\\
 & +\nabla\cdot\left(\nabla\left(u_{k}^{-a}v_{k}^{-b}\right)+\chi au_{k}^{-a}v_{k}^{-b-1}\nabla v_{k}\right)+b u_{k}^{-a}v_{k}^{-b}-b \frac{u_{k}^{{-a}+1}}{1+\frac{1}{k} u_{k}}v_{k}^{{-b}-1}.\label{KEYab}
\end{align}
Multiplying \eqref{KEYab} by an arbitrary function $\psi\in C^1_0(\R_0^+)$ and integrating by parts w.r.t. $t$ yields for all $x\in\Omega$ that
\begin{align}
&-\int_0^{\infty}u_{k}^{-a}v_{k}^{-b}\partial_t\psi\,dt-u_{k0}^{-a}v_{k0}^{-b}\psi(0)\nonumber\\
  = &-4\int_0^{\infty}\left(\frac{a+1}{a}\left|\nabla\left(u_{k}^{-\frac{{a}}{2}}v_{k}^{-\frac{{b}}{2}} \right)\right|^2+\left(\frac{b}{a}+\chi\frac{a+1}{2}\right)\nabla\left(u_{k}^{-\frac{{a}}{2}}v_{k}^{-\frac{{b}}{2}} \right)\cdot u_{k}^{-\frac{{a}}{2}}v_{k}^{-\frac{{b}}{2}-1}\nabla v_{k}\right.\nonumber\\
 &\phantom{-4aa\int_0^{\infty}}\left.+\frac{1}{4}\left(\frac{b^2}{a}+b+\chi b\right)\left|u_{k}^{-\frac{{a}}{2}}v_{k}^{-\frac{{b}}{2}-1}\nabla v_{k}\right|^2\right)\psi\,dt\nonumber\\
 & +\nabla\cdot\int_0^{\infty}\left(\nabla\left(u_{k}^{-a}v_{k}^{-b}\right)+\chi au_{k}^{-a}v_{k}^{-b-1}\nabla v_{k}\right)\psi\,dt+\int_0^{\infty}\left(b u_{k}^{-a}v_{k}^{-b}-b \frac{u_{k}^{{-a}+1}}{1+\frac{1}{k} u_{k}}v_{k}^{{-b}-1}\right)\psi.\label{KEYabt}
\end{align}
On the other hand, multiplying \eqref{KEYab} by an arbitrary function $\varphi\in C^1(\overline{\Omega})$ and integrating by parts w.r.t. $x$ and using the boundary conditions where necessary yields that 
\begin{align}
 \int_{\Omega}\partial_t\left(u_{k}^{-a}v_{k}^{-b}\right)\varphi\,=  &-4\int_{\Omega}\left(\frac{a+1}{a}\left|\nabla\left(u_{k}^{-\frac{{a}}{2}}v_{k}^{-\frac{{b}}{2}} \right)\right|^2+\left(\frac{b}{a}+\chi\frac{a+1}{2}\right)\nabla\left(u_{k}^{-\frac{{a}}{2}}v_{k}^{-\frac{{b}}{2}} \right)\cdot u_{k}^{-\frac{{a}}{2}}v_{k}^{-\frac{{b}}{2}-1}\nabla v_{k}\right.\nonumber\\
 &\left.\phantom{\int_{\Omega}-4}+\frac{1}{4}\left(\frac{b^2}{a}+b+\chi b\right)\left|u_{k}^{-\frac{{a}}{2}}v_{k}^{-\frac{{b}}{2}-1}\nabla v_{k}\right|^2\right)\varphi\,ds\nonumber\\
 &-\int_{\Omega}\left(\nabla\left(u_{k}^{-a}v_{k}^{-b}\right)+\chi au_{k}^{-a}v_{k}^{-b-1}\nabla v_{k}\right)\cdot\nabla\varphi\,ds\nonumber\\
 &+\int_{\Omega}\left(b u_{k}^{-a}v_{k}^{-b}-b u_{k}^{-a+1}v_{k}^{{-b}-1}\right)\varphi\,ds.
 \label{KEYabx}
\end{align}
Finally, if we multiply \eqref{KEYab} by the product $\psi\varphi$ and integrate by parts w.r.t. $t$ and $x$, then we arrive at the following variational reformulation of \eqref{KEYab}:
\begin{align}
 &-\int_0^{\infty}\partial_t\psi\int_{\Omega}u_{k}^{-a}v_{k}^{-b}\varphi\,dxdt-\psi(0)\int_{\Omega}u_{k0}^{-a}v_{k0}^{-b}\varphi\,dx\nonumber\\
  =  &-4\int_0^{\infty}\psi\int_{\Omega}\left(\frac{a+1}{a}\left|\nabla\left(u_{k}^{-\frac{{a}}{2}}v_{k}^{-\frac{{b}}{2}} \right)\right|^2+\left(\frac{b}{a}+\chi\frac{a+1}{2}\right)\nabla\left(u_{k}^{-\frac{{a}}{2}}v_{k}^{-\frac{{b}}{2}} \right)\cdot u_{k}^{-\frac{{a}}{2}}v_{k}^{-\frac{{b}}{2}-1}\nabla v_{k}\right.\nonumber\\
 &\left.\phantom{\int_0^t\psi\int_{\Omega}-4}+\frac{1}{4}\left(\frac{b^2}{a}+b+\chi b\right)\left|u_{k}^{-\frac{{a}}{2}}v_{k}^{-\frac{{b}}{2}-1}\nabla v_{k}\right|^2\right)\varphi\,dxds\nonumber\\
 &-\int_0^{\infty}\psi\int_{\Omega}\left(\nabla\left(u_{k}^{-a}v_{k}^{-b}\right)+\chi au_{k}^{-a}v_{k}^{-b-1}\nabla v_{k}\right)\cdot\nabla\varphi\,dxds\nonumber\\
 &+\int_0^{\infty}\psi\int_{\Omega}\left(b u_{k}^{-a}v_{k}^{-b}-b u_{k}^{-a+1}v_{k}^{{-b}-1}\right)\varphi\,dxds.
 \label{superue}
\end{align}
For equation \eqref{ve} a standard procedure yields the following reformulations:
for all $0\leq\varphi\in C^1(\overline{\Omega})$ and $0\leq\psi\in C^1_0(\R_0^+)$ it holds that
\begin{align}
 -\int_0^{\infty}v_{k}\partial_t\psi \,ds-\psi(0)v_{k0}=\nabla\cdot\int_0^{\infty}\psi\nabla v_{k}\,dt+\int_0^{\infty}\left(-v_{k}+\frac{u_{k}}{1+\frac{1}{k} u_{k}}\right)\psi\,dt.\label{supervte}
\end{align}
and
\begin{align}
 -\int_0^{\infty}\partial_t\psi\int_{\Omega}v_{k}\varphi\,dxds-\psi(0)\int_{\Omega}v_{k0}\varphi\,dx=\int_0^{\infty}\psi\int_{\Omega}-\nabla v_{k}\cdot\nabla\varphi+\left(-v_{k}+\frac{u_{k}}{1+\frac{1}{k} u_{k}}\right)\varphi\,dxds.\label{superve}
\end{align}
\subsection{Further uniform estimates for \texorpdfstring{\eqref{KSLoge}}{}}\label{keyest}
Choosing $\varphi\equiv 1$ in \eqref{KEYabx}  yields 
that 
\begin{align}
 \frac{d}{dt}\int_{\Omega} u_{k}^{-a}v_{k}^{-b}\,dx
 =&\int_{\Omega}-4\left(\frac{a+1}{a}\left|\nabla\left(u_{k}^{-\frac{{a}}{2}}v_{k}^{-\frac{{b}}{2}} \right)\right|^2+\left(\frac{b}{a}+\chi\frac{a+1}{2}\right)\nabla\left(u_{k}^{-\frac{{a}}{2}}v_{k}^{-\frac{{b}}{2}} \right)\cdot u_{k}^{-\frac{{a}}{2}}v_{k}^{-\frac{{b}}{2}-1}\nabla v_{k}\right.\nonumber\\
 &\phantom{\int_{\Omega}-4aa}\left.+\frac{1}{4}\left(\frac{b^2}{a}+b+\chi b\right)\left|u_{k}^{-\frac{{a}}{2}}v_{k}^{-\frac{{b}}{2}-1}\nabla v_{k}\right|^2\right)
 +b u_{k}^{-a}v_{k}^{-b}-b\frac{u_{k}^{{-a}+1}}{1+\frac{1}{k} u_{k}} v_{k}^{{-b}-1}\,dx.\label{Key}
\end{align}
Due to the key assumption \eqref{ab} it holds that 
\begin{align}
 &0>\left(\frac{b}{a}+\chi\frac{a+1}{2}\right)^2-\frac{a+1}{a}\left(\frac{b^2}{a}+b+\chi b\right)
 =\chi^2\frac{(a+1)^2}{4}-\frac{(b+a+1)b}{a}\nonumber\\
 \Leftrightarrow &\chi^2<4 \frac{(b+a+1)b}{a(a+1)^2}.\label{keyA}
\end{align}
Thus, the quadratic form
\begin{align}
 Q(U,V):=4\left(\frac{a+1}{a}\left|U\right|^2+\left(\frac{b}{a}+\chi\frac{a+1}{2}\right)U\cdot V+\frac{1}{4}\left(\frac{b^2}{a}+b+\chi b\right)\left|V\right|^2\right)\label{Q}
\end{align}
is strictly convex and satisfies
\begin{align}
 Q(U,V)\geq \Cl{low}\left(|U|^2+|V|^2\right)\qquad\text{ for all }U,V\in\R^n.\label{coer}
\end{align}
Then, due to \eqref{coer}, we have with \eqref{Key} that
\begin{align}
 \frac{d}{dt}\int_{\Omega} u_{k}^{-a}v_{k}^{-b}\,dx+\int_{\Omega}\Cr{low}\left(\left|\nabla\left(u_{k}^{-\frac{{a}}{2}}v_{k}^{-\frac{{b}}{2}} \right)\right|^2+\left|u_{k}^{-\frac{{a}}{2}}v_{k}^{-\frac{{b}}{2}-1}\nabla v_{k}\right|^2\right)+
 b \frac{u_{k}^{{-a}+1}}{1+\frac{1}{k} u_{k}}v_{k}^{{-b}-1}\,dx
 \leq&b\int_{\Omega} u_{k}^{-a}v_{k}^{-b}\,dx.\label{KeyApr}
\end{align}
Applying Gronwall's lemma to the differential inequality \eqref{KeyApr} we arrive at 
\begin{align}
 &\int_{\Omega} u_{k}^{-a}v_{k}^{-b}\,dx+\int_0^t\int_{\Omega}\Cr{low}\left(\left|\nabla\left(u_{k}^{-\frac{{a}}{2}}v_{k}^{-\frac{{b}}{2}} \right)\right|^2+\left|u_{k}^{-\frac{{a}}{2}}v_{k}^{-\frac{{b}}{2}-1}\nabla v_{k}\right|^2\right)+ b\frac{u_{k}^{{-a}+1}}{1+\frac{1}{k} u_{k}}v_{k}^{{-b}-1}\,dxds\nonumber\\
 \leq &e^{bt}\int_{\Omega} u_{k 0}^{-a}v_{k 0}^{-b}\,dx\nonumber\\
 \leq &\Cl{C2}(T).\label{est1}
\end{align}
Integral inequality \eqref{est1} yields the following set of estimates:
\begin{align}
 &\left\|\nabla\left(u_{k}^{-\frac{{a}}{2}}v_{k}^{-\frac{{b}}{2}} \right)\right\|_{L^2((0,T)\times\Omega)}\leq \Cl{C3}(T),\label{estU}\\
 &\left\|u_{k}^{-\frac{{a}}{2}}v_{k}^{-\frac{{b}}{2}-1}\nabla v_{k}\right\|_{L^2((0,T)\times\Omega)}\leq \Cr{C3}(T),\label{estV}\\
 &\left\|u_{k}^{-\frac{{a}}{2}} v_{k}^{-\frac{{b}}{2}}\right\|_{L^{\infty}(0,T;L^2(\Omega))}\leq \Cr{C3}(T),\label{estpr1}\\
 &\left\|\frac{u_{k}^{{-a}+1}}{1+\frac{1}{k} u_{k}}v_{k}^{{-b}-1}\right\|_{L^1((0,T)\times\Omega)}\leq\Cr{C3}(T).\label{est2}
\end{align}
Due to a Sobolev-type inequality (see, e.g., \cite[Chapter II \S 3 (3.4)]{LSU}), estimates  \eqref{estU} and \eqref{estpr1} imply that
\begin{align}
\left\|u_{k}^{-a} v_{k}^{-b}\right\|_{L^{1+\frac{2}{n}}((0,T)\times\Omega)}
 \leq &\C(T).\label{estpr3}
\end{align}
Thanks to estimates \eqref{estU} and \eqref{estpr3} and the H\"older inequality  we obtain that
\begin{align}
 &\left\|\nabla\left(u_{k}^{-a}v_{k}^{-b}\right)\right\|_{L^{\frac{n+2}{n+1}}((0,T)\times\Omega)}\leq \C(T).\label{estlord1}
\end{align}
Similarly, estimates \eqref{estV} and \eqref{estpr3}  and the H\"older inequality imply that
\begin{align}
 &\left\|u_{k}^{-a}v_{k}^{-b-1}\nabla v_{k}\right\|_{L^{\frac{n+2}{n+1}}((0,T)\times\Omega)}\leq \C(T).\label{estlord2}
\end{align}
Combining \eqref{estlord1}-\eqref{estlord2} we conclude that
\begin{align}
 \left\|\nabla\left(u_{k}^{-a}v_{k}^{-b}\right)+\chi au_{k}^{-a}v_{k}^{-b-1}\nabla v_{k}\right\|_{L^{\frac{n+2}{n+1}}((0,T)\times\Omega)}\leq \C(T).\label{estn2n1}
\end{align}
Thanks to  \eqref{estU}-\eqref{est2} and  \eqref{estn2n1} we deduce from \eqref{KEYabx}  that  
\begin{align}
 \left\|\partial_t\left(u_{k}^{-a}v_{k}^{-b}\right)\right\|_{L^1(0,T;(W^{1,n+2}(\Omega))^*)}\leq \Cl{C8}.\label{estprt}
\end{align}
Further, due to \eqref{aproxiniuv_} and  \eqref{estU}-\eqref{est2} we deduce from \eqref{KEYabt} that for all $\psi\in C^1_0(\R_0^+)$ it holds that
\begin{align}
 \left\|\nabla\cdot\int_0^{\infty}\psi\left(\nabla\left(u_{k}^{-a}v_{k}^{-b}\right)+\chi au_{k}^{-a}v_{k}^{-b-1}\nabla v_{k}\right)\,dt\right\|_{L^1(\Omega)}\leq \C(\psi).\label{estintnabla_}
\end{align}
Similarly, using \eqref{aproxiniv_} and \eqref{massue}-\eqref{massve} we deduce from \eqref{supervte} that
\begin{align}
 \left\|\nabla\cdot\int_0^{\infty}\psi\nabla v_{k}\,dt\right\|_{L^1(\Omega)}\leq \C(\psi).\label{estintDelta_}
\end{align}

\begin{Remark}[Sign of $b$]
 Condition \eqref{keyA} would also be satisfied if 
 \begin{align}
  b<b_-(a):=-\frac{1+a}{2}\left(\sqrt{1+\chi^2a}+1\right).\nonumber
 \end{align}
However, in this case 
\begin{align*}
b<0,
\end{align*}
 which implies that  the sign of the last term on the left-hand side of \eqref{KeyApr} is positive. Consequently, the immediate control upon this term is lost. Similar to the analysis \cite{LankWink2017} one would then be forced to impose restrictions upon $\chi$ and the space dimension $n$ in order to be able to bind it by means of others terms.
\end{Remark}

\section{Convergence to a generalised supersolution of \texorpdfstring{\eqref{KSLog}}{}: proof of {\it Theorem \ref{TheoEx}}}\label{SecEx}
Throughout this sequel we assume that 
\begin{align}
 0\leq\varphi\in C^1(\overline{\Omega}),\qquad 0\leq\psi\in C^1_0(\R_0^+),\nonumber
\end{align}
and are arbitrary.

Based on the properties established in {\it Section \ref{secreg}} we conclude that  a subsequence $k_m$ and  measurable functions $u,v:[0,\infty)\times \overline{\Omega}\rightarrow\R^+$ exist, such that:\\
due to \eqref{vcomp},
\begin{alignat}{3}
 &v_{k_m}\underset{m\rightarrow\infty}{\rightarrow}v&&\qquad\text{in }L^r_{loc}(\R_0^+\times\overline{\Omega})\text{ and a.e. in }\R^+\times\Omega,\label{convv}\\
 &\nabla v_{k_m}\underset{m\rightarrow\infty}{\rightarrow}\nabla v&&\qquad\text{in }L^s_{loc}(\R_0^+\times\overline{\Omega})\text{ and a.e. in }\R^+\times\Omega;\label{convnv}
\end{alignat}
due to \eqref{estvmin} and \eqref{convv}, for all $0<\tau<T<\infty$
\begin{align}
 \underset{(\tau,T)\times \Omega}{\ess\inf}v\geq \Cr{vmin}(\tau,T)>0,\label{vmin}
\end{align}
so that $v^{-1}\in L_{loc}^{\infty}(\R^+\times\overline{\Omega})$;\\
due to \eqref{estlord1}, \eqref{estprt}, and the Lions-Aubin lemma \cite[Corollary 4]{Simon},
\begin{align}
 u_{k_m}^{-a}v_{k_m}^{-b}\underset{m\rightarrow\infty}{\rightarrow}\eta\qquad\text{in }L^1_{loc}(\R_0^+\times\overline{\Omega})\text{ and a.e. in }\R^+\times\Omega;\label{convpr1_}
\end{align}
due to \eqref{convv} and \eqref{convpr1_}, 
\begin{align}
  u_{k_m}^{-1}\phantom{\ }=\phantom{\    }&\left(u_{k_m}^{-a}v_{k_m}^{-b}\right)^{\frac{1}{{a}}}v_{k_m}^{\frac{{b}}{{a}}}\nonumber\\
  \underset{m\rightarrow\infty}{\rightarrow}&\eta^{\frac{1}{{a}}}v^{\frac{{b}}{{a}}}\in[0,\infty)
  \qquad\text{  a.e. in }\R^+\times\Omega;\label{convum_}
\end{align}
due to \eqref{convum_},
\begin{align}
  u_{k_m}
  \underset{m\rightarrow\infty}{\rightarrow}&\eta^{-\frac{1}{{a}}}v^{-\frac{{b}}{{a}}}
  =:u\in(0,\infty]\qquad\text{  a.e. in }\R^+\times\Omega;\label{convu2}
\end{align}
due to \eqref{massue}, \eqref{convu2}, and Fatou's lemma,
\begin{align}
 &u_{k_m}
  \underset{m\rightarrow\infty}{\rightarrow}u\in (0,\infty)\qquad\text{  a.e. in }\R^+\times\Omega,\label{convu}
\end{align}
and $u\in L^{\infty}(\R^+;L^1(\Omega))$;\\
due to \eqref{convv}, \eqref{convpr1_},  \eqref{convu}, 
\begin{align}
 u_{k_m}^{-a}v_{k_m}^{-b}\underset{m\rightarrow\infty}{\rightarrow}u^{-a}v^{-b}\qquad\text{in }L^1_{loc}(\R_0^+\times\overline{\Omega})\text{ and a.e. in }\R^+\times\Omega;\label{convpr1}
\end{align}
due to \eqref{estV},  \eqref{convv}, \eqref{convnv},  \eqref{convu}, and the Lions lemma,
\begin{align}
 u_{k_m}^{-\frac{{a}}{2}}  v_{k_m}^{-\frac{{b}}{2}-1}\nabla v_{k_m}\underset{m\rightarrow\infty}{\rightharpoonup}u^{-\frac{{a}}{2}}  v^{-\frac{{b}}{2}-1}\nabla v\qquad\text{ in }L^2_{loc}(\R_0^+\times\overline{\Omega});\label{convwpr2}
\end{align}
due to \eqref{estU}, \eqref{convpr1},  and the Banach-Alaoglu theorem,
\begin{align}
 \nabla\left(u_{k_m}^{-\frac{{a}}{2}}v_{k_m}^{-\frac{{b}}{2}} \right)\underset{m\rightarrow\infty}{\rightharpoonup}\nabla\left(u^{-\frac{{a}}{2}}v^{-\frac{{b}}{2}}\right)\qquad\text{ in }L^2_{loc}(\R_0^+\times\overline{\Omega});\label{convwpr1}
\end{align}
due to \eqref{estlord1}, \eqref{convpr1},  and the Banach-Alaoglu theorem,
\begin{align}
 \nabla\left(u_{k_m}^{-a}v_{k_m}^{-b} \right)\underset{m\rightarrow\infty}{\rightharpoonup}\nabla\left(u^{-a}v^{-b}\right)\qquad\text{ in }L^{\frac{n+2}{n+1}}_{loc}(\R_0^+\times\overline{\Omega});\label{convwpr1_}
\end{align}
due to \eqref{estlord2}, \eqref{convv}, \eqref{convnv}, \eqref{convu},   and the Lions lemma,
\begin{align}
 u_{k_m}^{-a}v_{k_m}^{-b-1} \nabla v_{k_m}\underset{m\rightarrow\infty}{\rightharpoonup}u^{-a}v^{-b-1} \nabla v\qquad\text{ in }L^{\frac{n+2}{n+1}}_{loc}(\R_0^+\times\overline{\Omega});\label{convwpr1_2}
\end{align}
due to \eqref{convwpr1_}-\eqref{convwpr1_2},
\begin{align}
 \nabla\left(u_{k_m}^{-a}v_{k_m}^{-b} \right)+u_{k_m}^{-a}v_{k_m}^{-b-1} \nabla v_{k_m}\underset{m\rightarrow\infty}{\rightharpoonup}\nabla\left(u^{-a}v^{-b}\right)+u^{-a}v^{-b-1} \nabla v\qquad\text{ in }L^{\frac{n+2}{n+1}}_{loc}(\R_0^+\times\overline{\Omega});\label{conterm}
\end{align}
due to \eqref{convv} and \eqref{convu}, 
\begin{align}
 \frac{u_{k_m}^{{-a}+1}}{1+\frac{1}{k_m} u_{k_m}}v_{k_m}^{{-b}-1}
 \underset{m\rightarrow\infty}{\rightarrow}& u^{-a+1}v^{-b-1}\qquad\text{  a.e. in }\R^+\times\Omega;\label{conv4}
\end{align}
due to \eqref{est2}, \eqref{conv4}, and Fatou's lemma, 
\begin{align}
\underset{m\rightarrow\infty}{\lim\inf}\ \int_0^{\infty}\psi\int_{\Omega}\frac{u_{k_m}^{{-a}+1}}{1+\frac{1}{k_m} u_{k_m}}v_{k_m}^{{-b}-1}\varphi\,dxds
\geq&\int_0^{\infty}\psi\int_{\Omega} u^{-a+1}v^{{-b}-1}\varphi\,dxds,\label{convprodt}
\end{align}
and $u^{-a+1}v^{{-b}-1}\in L^1_{loc}(\R_0^+;L^1(\Omega))$;\\
due to  \eqref{convu}, 
\begin{align}
 \frac{u_{k_m}}{1+\frac{1}{k_m} u_{k_m}}
 \underset{m\rightarrow\infty}{\rightarrow}& u\qquad\text{  a.e. in }(0,T)\times\Omega;\label{conv5}
\end{align}
due to \eqref{aproxiniu_}, \eqref{massue}, and the Banach-Alaoglu theorem, there exists some $\mu\in {\cal M}(\R_0^+\times\overline{\Omega})$ s.t. 
\begin{align}
 \underset{m\rightarrow\infty}{\lim}\int_0^{\infty}\psi\int_{\Omega} \frac{u_{k_m}}{1+\frac{1}{k_m} u_{k_m}}\varphi\,dxds
 =\int_0^{\infty}\psi\int_{\Omega} u\varphi\,dxds+\int_0^{\infty}\psi\int_{\overline{\Omega}}\varphi\ d\mu(t,x);\label{conv6}
\end{align}
due to \eqref{conv5} and Fatou's lemma, 
\begin{align}
 \underset{m\rightarrow\infty}{\lim}\int_0^{\infty}\psi\int_{\Omega} \frac{u_{k_m}}{1+\frac{1}{k_m} u_{k_m}}\varphi\,dxds
 \geq\int_0^{\infty}\psi\int_{\Omega}u\varphi\,dxds;\label{conv7}
\end{align}
due to \eqref{conv6}  and  \eqref{conv7}, 
\begin{align}
 \int_0^{\infty}\psi\int_{\overline{\Omega}}\varphi\ d\mu(t,x)\geq0,\nonumber
\end{align}
so that $\mu$ is a non-negative measure;\\
due to \eqref{aproxiniu_} and  \eqref{massue},
\begin{align}
 \int_0^{\infty}\psi\int_{\Omega} \frac{u_{k_m}}{1+\frac{1}{k_m} u_{k_m}}\,dxds\phantom{\ }\leq\phantom{\    }&\int_0^{\infty}\psi\,ds\int_{\Omega}u_{k_m}\,dx\nonumber\\
 \phantom{\ }=\phantom{\    }&\int_0^{\infty}\psi\,ds\int_{\Omega}u_{k_m0}\,dx\nonumber\\
  \underset{m\rightarrow\infty}{\rightarrow}&\int_0^{\infty}\psi\,ds\int_{\Omega}u_0\,dx;\label{conv9}
\end{align}
due to \eqref{conv6} (set $\varphi\equiv 1$) and \eqref{conv9},  $u$ satisfies  \eqref{massu0} follows;\\
due to \eqref{convwpr2}, \eqref{convwpr1}, and  the convexity of $Q$ (was defined in \eqref{Q}), for all $0\leq\varphi\in C([0,T]\times\overline{\Omega})$
\begin{align}
 \underset{m\rightarrow\infty}{\lim\inf}\ &\int_0^T\int_{\Omega}\varphi Q\left(\nabla\left(u_{k_m}^{-\frac{{a}}{2}}v_{k_m}^{-\frac{{b}}{2}} \right),u_{k_m}^{-\frac{{a}}{2}}  v_{k_m}^{-\frac{{b}}{2}-1}\nabla v_{k_m}\right)\,dxds\nonumber\\
 \geq &\int_0^T\int_{\Omega}\varphi Q\left(\nabla\left(u^{-\frac{{a}}{2}}v^{-\frac{{b}}{2}} \right),u^{-\frac{{a}}{2}}  v^{-\frac{{b}}{2}-1}\nabla v\right)\,dxds;\label{convQ}
\end{align}
due to \eqref{estintnabla_},  \eqref{conterm}, and the Banach-Alaoglu theorem, 
\begin{align}
 &\nabla\cdot\int_0^{\infty}\psi\left(\nabla\left(u_{k_m}^{-a}v_{k_m}^{-b}\right)+\chi au_{k_m}^{-a}v_{k_m}^{-b-1}\nabla v_{k_m}\right)\,dt\nonumber\\
 \underset{m\rightarrow\infty}{\overset{*}{\rightharpoonup}}&\nabla\cdot\int_0^{\infty}\psi\left(\nabla\left(u^{-a}v^{-b}\right)+\chi au^{-a}v^{-b-1}\nabla v\right)\,dt\qquad\text{in }{\cal M}(\overline{\Omega}),\label{convbcpr}
\end{align}
and 
\begin{align*}
\int_0^{\infty}\psi\left(\nabla\left(u^{-a}v^{-b}\right)+\chi au^{-a}v^{-b-1}\nabla v\right)\,ds\in {\cal DM}^{\frac{n+2}{n+1}}(\Omega);   
\end{align*}
due to \eqref{contitr}, \eqref{uvebc}, \eqref{conterm}, and \eqref{convbcpr}, $(u,v)$ satisfies  \eqref{prweakbc};\\
due to \eqref{estintDelta_}, \eqref{convnv}, and the Banach-Alaoglu theorem, 
\begin{align}
 \nabla\cdot\int_0^{\infty}\psi \nabla v_{k_m}\,dt\underset{m\rightarrow\infty}{\overset{*}{\rightharpoonup}}\nabla\cdot\int_0^{\infty}\psi\nabla v\,dt\qquad\text{in }{\cal M}(\overline{\Omega}),\label{convDelta}
\end{align}
and 
\begin{align*}
\int_0^{\infty}\psi \nabla v\,ds\in {\cal DM}^{s}(\Omega);                                                              \end{align*}
due to \eqref{contitr}, \eqref{uvebc}, \eqref{convnv}, and \eqref{convDelta},  $v$ satisfies  \eqref{vweakbc}.

Altogether, combining \eqref{aproxiniuv_}, \eqref{convpr1},  \eqref{convprodt}, \eqref{convQ}, and  \eqref{convbcpr}, we deduce from \eqref{superue} for $k=k_m$ by taking limit superior as $m\rightarrow\infty$ that $(u,v)$ satisfies \eqref{superu}. Finally, passing to the limit in \eqref{superve} for $k=k_m$ as $m\rightarrow\infty$, and using \eqref{convv}, \eqref{convnv}, and \eqref{conv6}, we find that $(u,v)$ satisfies \eqref{superv}. {\it Theorem \ref{TheoEx}} is proved. 
\section{Classical solutions to \texorpdfstring{\eqref{KSLog}}{}: proof of {\it Theorem \ref{classol}}}\label{SecCl}
In this final section we prove  {\it Theorem \ref{classol}}. 
Thus, we assume now that  $u$ and $v$ are smooth and satisfy \eqref{udivv} and verify that in this case $(u,v)$ is, in fact, a classical solution to \eqref{KSLog}.

Since $v$ is smooth, the weak boundary condition \eqref{vweakbc} implies that $\partial_{\nu}v(t,\cdot)$ vanishes on the boundary of $\partial\Omega$ for all $t>0$. Further, since both $u$ and $v$ are smooth,  the weak boundary condition \eqref{prweakbc} takes the form
\begin{align}
 \partial_{\nu}\left(u^{-a}v^{-b}\right)+\chi au^{-a}v^{-b-1}\partial_{\nu} v=0\qquad\text{ for all }(t,x)\in(0,T)\times \partial\Omega.\nonumber
\end{align} 
This means that
\begin{align}
-au^{-a-1}v^{-b}\left(\partial_{\nu} u-\chi \frac{u}{v}\partial_{\nu} v\right)-bu^{-a}v^{-b-1}\partial_{\nu}v=0\qquad\text{ for all }(t,x)\in(0,T)\times \partial\Omega.\label{prclasbc}
\end{align}
Dividing  \eqref{prclasbc} by $-au^{-a-1}v^{-b}$ and plugging the boundary condition for $v$, we conclude that $\partial_{\nu}u(t,\cdot)$ vanishes on $\partial\Omega$ for all $t>0$ as well. 

Next, exploiting the smoothness of $u$ and $v$, we integrate by parts in  \eqref{superu} w.r.t. $t$ and $x$ and then apply the Du Bois-Reymond lemma. This results in the differential inequality 
\begin{align}
 \left(-au^{-a-1}v^{-b}\left(\partial_tu -\nabla\cdot\left(\nabla u-\chi \frac{u}{v}\nabla v\right)\right)\right)-bu^{-a}v^{-b-1}\left(\partial_t v-(\Delta v-v+u)\right)\leq0\qquad\text{in }\R^+\times\Omega.\label{Eq1ac}
\end{align} 
Similarly, \eqref{superv} implies that $\mu=\xi\,dxdt$ for some smooth density function $\xi\geq0$ and
\begin{align}
 \partial_t v=\Delta v-v+u+\xi\qquad\text{in }\R^+\times\Omega.\label{vxi}
\end{align}
Dividing \eqref{Eq1ac} by $-au^{-a-1}v^{-b}$ and making use of \eqref{vxi} we deduce that
\begin{align}
 &\partial_tu -\nabla\cdot\left(\nabla u-\chi \frac{u}{v}\nabla v\right)\geq-\frac{b}{a}\frac{u}{v}\xi\qquad\text{in }\R^+\times\Omega.\label{Eq1ineq}
 \end{align}
 Integrating \eqref{Eq1ineq} by parts over $\Omega$ thereby using the boundary conditions and, subsequently, integrating over $(0,t)$ for any $t>0$ we then have that
\begin{align}
 \|u(t,\cdot)\|_{L^1(\Omega)}\geq& \|u_0\|_{L^1(\Omega)}-\frac{b}{a}\int_0^t\int_{\Omega}\frac{u}{v}\xi\,dxds.\label{mu1}
\end{align}
On the other hand, due to the Du Bois-Reymond lemma inequality \eqref{massu0} is equivalent to the following: 
\begin{align}
 \|u(t,\cdot)\|_{L^1(\Omega)}+\|\xi(t,\cdot)\|_{L^1(\Omega)}\leq\|u_0\|_{L^1(\Omega)}\qquad\text{ for all }t\geq0.\label{mu2}
\end{align}
Thus, combining \eqref{mu1}-\eqref{mu2} with \eqref{udivv} we conclude that
\begin{align}
 \|\xi(t,\cdot)\|_{L^1(\Omega)}\leq&\frac{b}{a}\int_0^t\int_{\Omega}\frac{u}{v}\xi\,dxds\nonumber\\
 \leq& \C(T)\int_0^t\|\xi(s,\cdot)\|_{L^1(\Omega)}ds\qquad\text{ for all }0\leq t\leq T<\infty\nonumber
\end{align}
which yields that 
\begin{align}
\xi\equiv0\label{xi0} 
\end{align}
due to the Gronwall lemma. Plugging \eqref{xi0} into \eqref{vxi},  \eqref{Eq1ineq}, and \eqref{mu2} immediately yields that \eqref{Eqv} is satisfied in the classical sense, and it holds that
\begin{align}
 &\partial_tu -\nabla\cdot\left(\nabla u-\chi \frac{u}{v}\nabla v\right)\geq0\qquad\text{in }\R^+\times\Omega\label{Eq1ineq0}
 \end{align}
 and
 \begin{align}
  \|u(t,\cdot)\|_{L^1(\Omega)}=& \|u_0\|_{L^1(\Omega)}\qquad\text{for all }t\geq0.\label{mas}
 \end{align}
Since \eqref{Eq1ineq0} is subject to the no-flux boundary conditions,   \eqref{Eq1ineq0}  and \eqref{mas} imply that equality holds in \eqref{Eq1ineq0}. Thus, equation \eqref{Equ} is satisfied in the classical sense, and {\it Theorem \ref{classol}} is proved.
\begin{appendices}
\section{Divergence-measure fields and their normal traces}\label{AppA}
In this section we collect some facts concerning the  divergence-measure fields.


Let ${\cal M}(\overline{\Omega})$  denote the space of Radon measures in $\overline{\Omega}$.  We recall the definition of the Banach space of divergence-measure fields \cite{ChenFrid2001} and its norm: 
\begin{align}
&{\cal DM}^p(\Omega):=\{F\in (L^p(\Omega))^n|\ \ \nabla\cdot F\in {\cal M}(\overline{\Omega})\},\nonumber\\
 &\|F\|_{{\cal DM}^p(\Omega)}:=\|F\|_{(L^p(\Omega))^n}+\|\nabla\cdot F\|_{{\cal M}(\overline{\Omega})}.\nonumber
\end{align}
Thereby we assume that $p\in\left(1,\frac{n}{n-1}\right)$, which is sufficient for our needs. 
Following \cite{ChenFrid2001} we introduce a generalisation of the normal trace over the boundary of $\partial\Omega$ which automatically satisfies a  Gauss-Green formula:  
\begin{align}
 \left<F\cdot\nu|_{\partial\Omega},\varphi\right>:=\int_{\Omega}F\cdot\nabla ({\cal E}\varphi)\,dx+\int_{\overline{\Omega}} ({\cal E}\varphi)\, d(\nabla\cdot F)\qquad\text{for all }\varphi\in  W^{\frac{1}{p},\frac{p}{p-1}}(\partial\Omega).\label{GG}
\end{align}
Here ${\cal E}: W^{\frac{1}{p},\frac{p}{p-1}}(\partial\Omega)\rightarrow W^{1,\frac{p}{p-1}}(\Omega)$ is a usual   extension operator, i.e., a continuous right inverse of the  corresponding trace operator. 
  It is known (see \cite[Theorem 2.1]{ChenFrid2001}) that  $F\cdot\nu|_{\partial\Omega}\in W^{-\frac{1}{p},p}(\partial\Omega)$ and  doesn't depend upon the particular choice of ${\cal E}$. Formula \eqref{GG} ensures the following implication:
  \begin{align}
   \begin{cases}F_n\underset{m\rightarrow\infty}{{\rightharpoonup}}F&\text{ in }(L^p(\Omega))^n,\\ \nabla \cdot F_n\underset{m\rightarrow\infty}{\overset{*}{\rightharpoonup}}\nabla\cdot F&\text{ in }{\cal M}(\overline{\Omega})
   \end{cases}\quad\Rightarrow\quad  F_n\cdot\nu|_{\partial\Omega}\underset{m\rightarrow\infty}{\overset{*}{\rightharpoonup}}F\cdot\nu|_{\partial\Omega}\qquad\text{in }W^{-\frac{1}{p},p}(\partial\Omega).\label{contitr}
  \end{align}

%
\end{appendices}
\phantomsection
\printbibliography
\end{document}